\documentclass{article}
\usepackage{amsmath,amssymb,amsthm}
\usepackage{parskip}
\usepackage{authblk}
\usepackage[numbers]{natbib}
\usepackage{geometry}
\usepackage{multicol}
\usepackage{enumerate}
\usepackage[bookmarks=false,colorlinks=true,citecolor=blue]{hyperref}
\usepackage{tikz}
\usepackage{orcidlink}
\usetikzlibrary{automata,positioning}

\newcommand\setItemnumber[1]{\setcounter{enumi}{\numexpr#1-1\relax}}

\newtheorem{theorem}{Theorem}[section]{\bfseries}{\itshape}
\newtheorem{definition}[theorem]{Definition}{\bfseries}{\itshape}
\newtheorem{corollary}[theorem]{Corollary}{\bfseries}{\itshape}
\newtheorem{puzzle}[theorem]{Puzzle}{\bfseries}{\itshape}
\newtheorem{case}{Case}
\newtheorem{subcase}{Subcase}
\numberwithin{subcase}{case}

\numberwithin{subsubcase}{subcase}

\title{\textbf{A Graph-Theoretic Model for a Generic Three-Jug Puzzle}}
\author[1]{\textbf{Suresh~Manjanath~Hegde}\orcidlink{0000-0003-0788-784X}}
\author[1,2]{\textbf{Shashanka~Kulamarva}\orcidlink{0009-0002-2982-6044}}
\affil[1]{National Institute of Technology Karnataka, Surathkal-575025, India}
\affil[2]{Graduate School of Informatics, Kyoto University, Kyoto-6068501, Japan}
\affil[ ]{Email: \textit{smhegde@nitk.edu.in, kulamarva.shashanka.3k@kyoto-u.ac.jp}}
\date{}

\begin{document}
	\maketitle
	\begin{abstract}
		\noindent A classic three-jug puzzle asks, given three jugs $A$, $B$, and $C$ with fixed maximum capacities, with jug $A$ filled with wine to its maximum capacity, whether is it possible to divide the wine into two halves by pouring it from one jug to another without using any other measuring devices. However, we consider a generic version of the three-jug puzzle and present an independent graph-theoretic model to determine whether the puzzle has a solution at all. If it has a solution, then the same can be determined using this model. We also present the sketch of an algorithm to determine the solution of the puzzle.\\
		
		
		
	\end{abstract}
	
	\section{Introduction}
	Mathematics can be a powerful tool for solving various puzzles, riddles, and open questions. In most cases, mathematical modeling helps in this regard. A three-jug puzzle is one of the ancient and familiar puzzles that derived out of pure recreational curiosity but later resulted in significant computational techniques. A formal statement of the puzzle goes as follows:
	
	\begin{puzzle}[\textbf{Three-Jug Puzzle}]
		Given three jugs $A$, $B$, and $C$ having the maximum integer capacities $a$, $b$, and $c$ gallons respectively (hereafter, we use the word `capacity' in place of `maximum integer capacity') with $a>b>c$, and the jug $A$ being filled with wine to its capacity with jugs $B$ and $C$ being empty, is it possible to divide the wine into two halves only by pouring it from one jug to another without using any other measuring devices? An instance of a three-jug puzzle is denoted by $\mathbb{P}_{a,b,c}$.
	\end{puzzle}
	
	Alternatively, the puzzle can be framed as the following question: A milkman has a jug (of capacity $a$ gallons with $a$ being an even number) full of milk. The customer has two jugs of capacities $b$ and $c$. The milkman wants to provide the customer with exactly $\frac{a}{2}$ gallons of milk. Is this possible even though the milkman does not have any measuring instruments? For a better understanding of the puzzle, we come up with a definition of a specific way of pouring that is in some sense `measurable'.
	
	\begin{definition}[\textbf{Measurable Pour}]
		A pour of wine from jug $X$ to jug $Y$ ($X$ and $Y$ can be any two jugs out of $A$, $B$, and $C$), is said to be a measurable pour, denoted by $X \downarrow Y$, if the pouring is stopped when the jug $X$ becomes empty or the jug $Y$ reaches its capacity.
	\end{definition}
	
	Note that a measurable pour keeps track of the volumes of wine in each jug after pouring; provided we know the volumes before pouring. One can see that if we don't perform a measurable pour and stop the pouring somewhere in between, we will loose the information on the volume of wine in each jug since we can not use any measuring devices. Since we need to keep track of the amount of wine in each jug without using any measuring devices, every pour should be a measurable pour.
	
	The earliest recorded appearance of jug-based measurement puzzles can be traced back to 1484 in \emph{Triparty en la science des nomvres}, a work of a famous French mathematician, Nicolas Chuquet \citep[see][]{Petkovic2000Puzzles}. Later, the puzzle was studied by Tartaglia, an Italian mathematician of the 16th century \cite[see][]{McDiarmid1994Jugs}. At the beginning of the 20th century, in Wilhelm Ahrens' foundational work in recreational mathematics, \emph{Mathematische Unterhaltungen und Spiele} \cite{Ahrens1901Games}, puzzles of constraint-based reasoning, including jug-based measurement problems, were framed as mathematical amusements.
	
	\citet{Tweedie1939MeasuringPuzzle} proposed a graph-theoretic solution for visualizing state transitions in Tartaglian puzzles, which later directly led to the current approaches to the three-jug problem \cite{Ore1990GraphsUses}. The puzzle came into the limelight again in the 1940s when \citet{Grossman1940GenPuzzle} extended the puzzle's formulation, characterizing a broader family of liquid-transfer problems. \citet{Buker1941Problems} formally introduced the Three-Jug Puzzle as \emph{Problem E451} in `Problems and Solutions', prompting a series of scholarly responses, including notable partial solutions by Browne \cite{Browne1942Problems} and \citet{Currie1946Problems}, and an insightful generalization by \citet{Scott1944Problems}. However, in these works, the focus shifted onto determining the values of $a,b,$ and $c$ (jug capacities) for which the solution exists.
	
	To begin with, researchers intuitively sought number-theoretic and algebraic models to formulate these puzzles. \citet{Goodstein1941MeasuringProblem} discussed the two-jug case (no jug C; pouring is only between the jugs A and B or emptying or filling one of them) presenting it as a Diophantine equation problem. \citet{Lawrence1941MeasuringProblem} presented an algebraic framework for analyzing pouring problems, emphasizing state-space exploration via modular arithmetic and linear combinations. \citet{McDiarmid1994Jugs} considered the three-jug puzzle with a condition $a=b+c$ and solved the puzzle using some combinatorial strategies wherein the jug contents were modeled within lattice frameworks, enabling solutions through integer partitioning techniques. They also proved that a solution to the puzzle exists if and only if the number $a$ is divisible by $2r$ where $r=gcd(b,c)$. However, the puzzle becomes more interesting without the condition $a=b+c$.
	
	There were geometric and computational approaches proposed for the puzzle. \citet{Lalchev2009GeomProblem} proposed Perlman’s geometric method for the puzzle, which visually represents possible jug states as vectors within a geometric space. \citet{Goetschalckx2011Flow} approached the jug problem from an operations research perspective, viewing it as a special case of flow routing in supply networks.
	
	The formal algorithmic approaches for the water-jug puzzles began when \citet{Boldi2002Jugs} defined the state-space graph and proved complexity bounds for the same together with a study on the solvability of the puzzle. \citet{Leon2005GenericTwoJugs} studied a generic two-jug puzzle and proposed a solution by designing heuristic algorithms, using artificial intelligence principles and state exploration strategies. Later, \citet{Man2015TwoJugAlgo} developed a rigorous algorithmic approach for solving the two-jug problem, which generalizes to higher dimensions. \citet{OBeirne1965Puzzles} and \citet{Ore1990GraphsUses} included jug puzzles as examples of mathematical paradoxes and graph-based challenges. \citet{Murray2003Puzzle} suggested that the puzzle could be approached computationally since it utilizes state-space representation and heuristic reasoning. \citet{Petkovic2000Puzzles}, in his book, compiled jug puzzles alongside other classical brainteasers, emphasizing their historical context and the mathematical creativity they evoke. For a historical review on the three-jug puzzle, one can refer to \cite{OBeirne1965Puzzles}.
	
	Inspired by Tweedie's approach \cite{Tweedie1939MeasuringPuzzle}, \citet{Ore1990GraphsUses} proposed a model for the three-jug puzzle in the framework of graph theory by using directed paths in directed graphs. Any distribution of the wine in the three jugs $A$, $B$, and $C$, can be described by the quantities $b$ and $c$ of wine in the jugs $B$ and $C$, respectively. Hence, every possible distribution of wine is denoted by a pair $(b,c)$.
	
	\begin{figure}[h]
		\centering
		\begin{tikzpicture}[scale=0.70]%
			[>=stealth,
			shorten >=1pt,
			node distance=2cm,
			on grid,
			auto,
			every state/.style={draw=black!60, fill=black!5, very thick}
			]
			\begin{scope}[every node/.style={circle,draw,inner sep=0pt}]
				\node (a) at (0,0) {$(0,0)$};
				\node (b) at (0,2) {$(0,1)$};
				\node (c) at (0,4) {$(0,2)$};
				\node (d) at (0,6) {$(0,3)$};
				\node (e) at (2,0) {$(1,0)$};
				\node (f) at (2,2) {$(1,1)$};
				\node (g) at (2,4) {$(1,2)$};
				\node (h) at (2,6) {$(1,3)$};
				\node (i) at (4,0) {$(2,0)$};
				\node (j) at (4,2) {$(2,1)$};
				\node (k) at (4,4) {$(2,2)$};
				\node (l) at (4,6) {$(2,3)$};
				\node (m) at (6,0) {$(3,0)$};
				\node (n) at (6,2) {$(3,1)$};
				\node (o) at (6,4) {$(3,2)$};
				\node (p) at (6,6) {$(3,3)$};
				\node (q) at (8,0) {$(4,0)$};
				\node (r) at (8,2) {$(4,1)$};
				\node (s) at (8,4) {$(4,2)$};
				\node (t) at (8,6) {$(4,3)$};
				\node (u) at (10,0) {$(5,0)$};
				\node (v) at (10,2) {$(5,1)$};
				\node (w) at (10,4) {$(5,2)$};
				\node (x) at (10,6) {$(5,3)$};
				\node (y) at (12,0) {$(6,0)$};
				\node (z) at (12,2) {$(6,1)$};
				\node (aa) at (12,4) {$(6,2)$};
				\node (ab) at (12,6) {$(6,3)$};
				\node (ac) at (14,0) {$(7,0)$};
				\node (ad) at (14,2) {$(7,1)$};
				\node (ae) at (14,4) {$(7,2)$};
				\node (af) at (14,6) {$(7,3)$};
			\end{scope}
			\draw[thick, ->, brown] (a) to (ac);
			\draw[thick, ->, green] (ac) to (t);
			\draw[thick, ->, blue] (t) to (q);
			\draw[thick, ->, yellow] (q) to (h);
			\draw[thick, ->, orange] (h) to (e);
			\draw[thick, ->, violet] (e) to (b);
			\draw[thick, ->, cyan] (b) to (ad);
			\draw[thick, ->, magenta] (ad) to (x);
			\draw[thick, ->, red] (x) to (u);
		\end{tikzpicture}
		\caption{Directed path indicating a solution to the three-jug puzzle $\mathbb{P}_{10,7,3}$}
		\label{fig:ThreeJug}
	\end{figure}
	
	Let us consider an illustration with the three-jug puzzle $\mathbb{P}_{10,7,3}$. Initially, all the wine is in jug $A$ indicating that $b=c=0$ so that the initial distribution is $(0,0)$. The target distribution to be achieved is $(5,0)$. The idea is to consider all non-negative integer pairs from $(0,0)$ to $(b,c)$ as vertices of a graph and to draw an edge wherever it's possible to move from one vertex (distribution) to another vertex (distribution) by a measurable pour of wine between a pair of jugs. This results in a directed path in the graph which indicates a solution to the three-jug puzzle as depicted in the Figure~\ref{fig:ThreeJug}. A directed $uv$-path means a directed path from a vertex $u$ to a vertex $v$.
	
	Observe that jug $A$ being filled with wine to its capacity while jugs $B$ and $C$ are empty is a restriction on the three-jug puzzle. The puzzle becomes more interesting when this condition is removed. Therefore, we remove this restriction and consider a generic version of the three-jug puzzle defined as follows.
	
	\begin{puzzle}[\textbf{Generic Three-Jug Puzzle}]
		We are given with three jugs $A$, $B$, and $C$ having the capacities $a$, $b$, and $c$ gallons respectively, with $a > b > c$. Currently, the three jugs are filled with $\tilde{a}$, $\tilde{b}$, and $\tilde{c}$ gallons of wine, respectively. Note that all the numbers are positive integers. Let $d:=\tilde{a}+\tilde{b}+\tilde{c}$, and we are given that for some positive integer $k$, $d=2k$ and $b \ge k$. Is it possible to divide the $d$ gallons of wine into two halves only by pouring (i.e., a series of measurable pours) it from one jug to another without using any other measuring devices? An instance of a generic three-jug puzzle is denoted by $\mathbb{P}^{\tilde{a},\tilde{b},\tilde{c}}_{a,b,c}$ (see Figure~\ref{fig:GenThreeJug}).
	\end{puzzle}
	
	\begin{figure}[h]
		\centering
		\begin{tikzpicture}[scale=0.70]
			
			% Jug dimensions and spacing
			\def\jugWidth{2}
			\def\jugGap{4}
			
			% Heights (capacities)
			\def\Ha{6}
			\def\Hb{4.5}
			\def\Hc{3.5}
			
			% Current fill heights
			\def\HaPrime{4.2}
			\def\HbPrime{3}
			\def\HcPrime{1.5}
			
			% Jug A
			\begin{scope}[xshift=0cm]
				\draw[thick] (0,0) rectangle (\jugWidth,\Ha);
				\fill[blue!30] (0,0) rectangle (\jugWidth,\HaPrime);
				\node at (\jugWidth/2, -0.7) {Jug A};
				\draw[dashed] (0,\Ha) -- (\jugWidth,\Ha);
				\node[right] at (\jugWidth, \Ha) {\(a\)};
				\node at (\jugWidth/2, \HaPrime + 0.3) {\(\tilde{a}\)};
			\end{scope}
			
			% Jug B
			\begin{scope}[xshift=\jugGap cm]
				\draw[thick] (0,0) rectangle (\jugWidth,\Hb);
				\fill[blue!30] (0,0) rectangle (\jugWidth,\HbPrime);
				\node at (\jugWidth/2, -0.7) {Jug B};
				\draw[dashed] (0,\Hb) -- (\jugWidth,\Hb);
				\node[right] at (\jugWidth, \Hb) {\(b\)};
				\node at (\jugWidth/2, \HbPrime + 0.3) {\(\tilde{b}\)};
			\end{scope}
			
			% Jug C
			\begin{scope}[xshift=2*\jugGap cm]
				\draw[thick] (0,0) rectangle (\jugWidth,\Hc);
				\fill[blue!30] (0,0) rectangle (\jugWidth,\HcPrime);
				\node at (\jugWidth/2, -0.7) {Jug C};
				\draw[dashed] (0,\Hc) -- (\jugWidth,\Hc);
				\node[right] at (\jugWidth, \Hc) {\(c\)};
				\node at (\jugWidth/2, \HcPrime + 0.3) {\(\tilde{c}\)};
			\end{scope}
			
		\end{tikzpicture}
		\caption{A generic three-jug puzzle $\mathbb{P}^{\tilde{a},\tilde{b},\tilde{c}}_{a,b,c}$}
		\label{fig:GenThreeJug}
	\end{figure}
	
	Observe that in the generic three-jug puzzle, the condition that for some positive integer $k$, $d=2k$ and $b \ge k$ are the minimum requirements for the existence of a solution. If any one of those does not hold, there is no solution for the puzzle.
	
	The goal of this paper is to solve this generic three-jug puzzle efficiently in a unique way. In this regard, inspired by \citet{Ore1990GraphsUses}, we propose a graph-theoretic solution to the generic three-jug puzzle by coming up with a directed graph model to check whether the given generic three-jug puzzle is solvable or not. If the puzzle has a solution, then it can be determined using the same model.
	
	\section{The Graph Model}\label{sec:Model}
	
	This section aims to come up with a directed graph for a given `special' ordered quadruple of positive integers. In this regard, we have the following definition.
	
	\begin{definition}\label{def:Model}
		For a given ordered quadruple $Q=(a,b,c,d)$ of positive integers that satisfy the following conditions:
		\begin{itemize}
			\item $a > b > c$
			\item For some positive integer $k$, $d=2k$ and $b \ge k$,
		\end{itemize}
		a corresponding directed graph $G_Q$ is defined as $G_Q=(V(G_Q),E(G_Q))$ with the set of vertices $V(G_Q)$ being the set of all ordered pairs $(i,j)$ with $0 \le i \le b$ and $0 \le j \le c$, i.e.,
		\begin{equation*}
			V(G_Q)=\bigcup_{\substack{0 \le i \le b \\ \mathstrut 0 \le j \le c \\ \mathstrut i,j \in \mathbb{Z}}}(i,j)
		\end{equation*}
		and there exists a directed edge in $E(G_Q)$ from a vertex $(i,j) \in V(G_Q)$ to another vertex $(i',j') \in V(G_Q)$ if and only if $0 \le d-i-j \le a$, $0 \le d-i'-j' \le a$, and one of the following conditions hold:
		\begin{enumerate}
			\item $j'=j$ and $i' \in \{0,b,d-j,d-a-j\}$
			\item $i'=i$ and $j' \in \{0,c,d-i,d-a-i\}$
			\item $i+j=i'+j'$ and \{$i' \in \{0,b\}$ or $j' \in \{0,c\}$\}.
		\end{enumerate}
	\end{definition}
	
	An illustration of graph $G_Q$ generated by the quadruple $Q = (7,4,2,6)$ is given in Figure~\ref{fig:ModelExample}. Note that the final graph obtained by our model looks similar to the graph obtained by \citet{Ore1990GraphsUses}; but there are a couple of significant differences implying the advantages of our graph model.
	
	First, this model is used to tackle a puzzle (reasoning in the next section) that is a more relaxed, generic version of the three-jug puzzle which was considered in \cite{Ore1990GraphsUses}. Further, in the graph model in \cite{Ore1990GraphsUses}, the existence of the edges depends on the decision whether one can achieve a distribution from another distribution by pouring wine between the jugs; whereas our graph modeling does not require the prior information of whether it is possible to achieve a target distribution from a source distribution by a valid pour or not. Moreover, this graph model is completely independent of the puzzle itself since we did not use the generic three-jug puzzle while defining this graph model and all that we need is an ordered quadruple satisfying some conditions. Moreover, our way of modeling nicely comes out as an algorithm.
	
	\begin{figure}[h]
		\centering
		\begin{tikzpicture}[scale=0.90]%
			[>=stealth,
			shorten >=1pt,
			node distance=2cm,
			on grid,
			auto,
			every state/.style={draw=black!60, fill=black!5, very thick}
			]
			\begin{scope}[every node/.style={circle,draw,inner sep=0pt}]
				\node (a) at (0,0) {$(0,0)$};
				\node (b) at (0,2) {$(0,1)$};
				\node (c) at (0,4) {$(0,2)$};
				\node (d) at (2,0) {$(1,0)$};
				\node (e) at (2,2) {$(1,1)$};
				\node (f) at (2,4) {$(1,2)$};
				\node (g) at (4,0) {$(2,0)$};
				\node (h) at (4,2) {$(2,1)$};
				\node (i) at (4,4) {$(2,2)$};
				\node (j) at (6,0) {$(3,0)$};
				\node (k) at (6,2) {$(3,1)$};
				\node (l) at (6,4) {$(3,2)$};
				\node (m) at (8,0) {$(4,0)$};
				\node (n) at (8,2) {$(4,1)$};
				\node (o) at (8,4) {$(4,2)$};
			\end{scope}
			\draw[<->, bend left] (a) to (c);
			\draw[<->, bend right] (a) to (m);
			\draw[->] (b) to (a);
			\draw[->] (b) to (c);
			\draw[<->, bend right] (b) to (n);
			\draw[<->, bend left] (c) to (o);
			\draw[->] (d) to (a);
			\draw[<->] (d) to (b);
			\draw[->, bend right] (d) to (m);
			\draw[<->, bend left] (d) to (f);
			\draw[->] (e) to (b);
			\draw[->] (e) to (c);
			\draw[->] (e) to (d);
			\draw[->] (e) to (f);
			\draw[->] (e) to (g);
			\draw[->, bend left] (e) to (n);
			\draw[->] (f) to (c);
			\draw[->, bend left] (f) to (o);
			\draw[->, bend left] (g) to (a);
			\draw[<->, bend left] (g) to (c);
			\draw[->, bend right] (g) to (m);
			\draw[<->, bend left] (g) to (i);
			\draw[->, bend left] (h) to (b);
			\draw[->] (h) to (f);
			\draw[->] (h) to (g);
			\draw[->] (h) to (i);
			\draw[->] (h) to (j);
			\draw[->, bend left] (h) to (n);
			\draw[->, bend right] (i) to (c);
			\draw[->, bend left] (i) to (o);
			\draw[->, bend left] (j) to (a);
			\draw[<->, bend left] (j) to (f);
			\draw[->] (j) to (m);
			\draw[<->, bend left] (j) to (l);
			\draw[->, bend left] (k) to (b);
			\draw[->] (k) to (i);
			\draw[->] (k) to (j);
			\draw[->] (k) to (l);
			\draw[->] (k) to (m);
			\draw[->] (k) to (n);
			\draw[->, bend right] (l) to (c);
			\draw[->] (l) to (o);
			\draw[<->, bend left] (m) to (i);
			\draw[<->, bend right] (m) to (o);
			\draw[<->] (n) to (l);
			\draw[->] (n) to (m);
			\draw[->] (n) to (o);
		\end{tikzpicture}
		\caption{Graph $G_Q$ corresponding to the quadruple $Q=(7,4,2,6)$}
		\label{fig:ModelExample}
	\end{figure}
	
	\section{Main Results}
	
	While defining the graph model in the previous section, the generic three-jug puzzle is nowhere in the picture. However, if we consider any generic three-jug puzzle $\mathbb{P}^{\tilde{a},\tilde{b},\tilde{c}}_{a,b,c}$, it indirectly generates an ordered quadruple $P=(a,b,c,\tilde{a}+\tilde{b}+\tilde{c})$; this ordered quadruple $P$ satisfies the preliminary requirements of the model which in turn results in the generation of a directed graph $G_P$ as mentioned in Definition~\ref{def:Model}. Note that there is a one-to-one correspondence between a vertex $(i,j)$ in $G_P$ and a distribution $(i,j)$ in the puzzle $\mathbb{P}^{\tilde{a},\tilde{b},\tilde{c}}_{a,b,c}$.
	
	Moving ahead, the idea is to utilize this graph $G_P$ and its properties to look for a solution to the generic three-jug puzzle $\mathbb{P}^{\tilde{a},\tilde{b},\tilde{c}}_{a,b,c}$, if it exists; but how can this be done? This section is all about answering this question. In order to do that, we need to prove a theorem and its corollary that are provided in the section.
	
	\begin{theorem}\label{thm:ThreeJug}
		Let $\mathbb{P}^{\tilde{a},\tilde{b},\tilde{c}}_{a,b,c}$ be a generic three-jug puzzle and let $G_P$ be the graph obtained from the ordered quadruple $P=(a,b,c,\tilde{a}+\tilde{b}+\tilde{c})$ using the graph model in Definition~\ref{def:Model}. Then there exists a directed edge from a vertex $(i,j)$ to a vertex $(i',j')$ in $G_P$ if and only if the corresponding distribution $(i',j')$ can be obtained from the distribution $(i,j)$ in the puzzle $\mathbb{P}^{\tilde{a},\tilde{b},\tilde{c}}_{a,b,c}$ by a single measurable pour.
	\end{theorem}
	
	\begin{proof}
		We have the generic three-jug puzzle $\mathbb{P}^{\tilde{a},\tilde{b},\tilde{c}}_{a,b,c}$, the ordered quadruple $P=(a,b,c,\tilde{a}+\tilde{b}+\tilde{c})$, and the graph $G_P$ as in the theorem statement. Let $d=\tilde{a}+\tilde{b}+\tilde{c}$. First, assume that the distribution $(i',j')$ can be obtained from the distribution $(i,j)$ in the puzzle $\mathbb{P}^{\tilde{a},\tilde{b},\tilde{c}}_{a,b,c}$ by a single measurable pour $X \downarrow Y$. Now, the goal is to prove that:
		\begin{equation}\tag{$\star$}
			\text{there exists a directed edge from the vertex } (i,j) \text{ to the vertex } (i',j') \text{ in } G_P\text{.}
		\end{equation}
		Since $d-i-j$ and $d-i'-j'$ are the levels of wines in the jug $A$ in the distribution $(i,j)$ and $(i',j')$ respectively, we have:
		\begin{equation*}
			0 \le d-i-j \le a \text{\hspace{0.75cm} and \hspace{0.75cm}} 0 \le d-i'-j' \le a
		\end{equation*}
		Depending on the choice of jugs $X$ and $Y$, we have the following cases of measurable pours possible.
		
		\begin{case}
			$X=A$ and $Y=B$, i.e., the measurable pour is $A \downarrow B$.
		\end{case}
		
		In this case, since the jug $C$ is untouched, we have $j'=j$. Since a measurable pour is being performed, we stop pouring only if jug $A$ becomes empty or jug $B$ becomes full. Depending on the reason for the termination of pouring, we have the following subcases.
		
		\begin{subcase}
			Pouring terminated because the jug $A$ became empty.
		\end{subcase}
		
		In this case, since jug $A$ is empty, $d$ gallons of wine is distributed only in jugs $B$ and $C$ which implies that $i'+j'=d$. This together with the fact that $j'=j$ implies that $i'=d-j$. Thus we have $i' \in \{0,b,d-j,d-a-j\}$.
		
		\begin{subcase}
			Pouring terminated because the jug $B$ became full.
		\end{subcase}
		
		In this case, since jug $B$ is full, the amount of wine in jug $B$ is equal to its capacity. Hence, we have that $i'=b$ indicating that $i' \in \{0,b,d-j,d-a-j\}$.
		
		In any subcase, we have $i' \in \{0,b,d-j,d-a-j\}$, collectively implying $(\star)$, as desired.
		
		\begin{case}
			The measurable pour is $B \downarrow A$.
		\end{case}
		
		In this case, since the jug $C$ is untouched, we have $j'=j$. Since a measurable pour is being performed, we stop pouring only if jug $B$ becomes empty or jug $A$ becomes full. Again, depending on the reason for the termination of pouring, we have the following subcases.
		
		\begin{subcase}
			Pouring terminated because the jug $B$ became empty.
		\end{subcase}
		
		In this case, since the jug $B$ is empty, the amount of wine in the jug $B$ is equal to $0$. Hence, we have $i'=0$ indicating that $i' \in \{0,b,d-j,d-a-j\}$.
		
		\begin{subcase}
			Pouring terminated because the jug $A$ became full.
		\end{subcase}
		
		In this case, since jug $A$ is full, the amount of wine in jug $A$ is equal to its capacity, and the remaining wine is distributed over jugs $B$ and $C$ which implies that $a+i'+j'=d$. Further, substituting $j'=j$, we get $a+i'+j=d$. Hence, we have that $i'=d-a-j$ indicating that $i' \in \{0,b,d-j,d-a-j\}$.
		
		In any subcase, we have $i' \in \{0,b,d-j,d-a-j\}$, collectively implying $(\star)$, as desired.
		
		\begin{case}
			The measurable pour is $A \downarrow C$.
		\end{case}
		
		In this case, since the jug $B$ is untouched, we have $i'=i$. Since a measurable pour is being performed, we stop pouring only if jug $A$ becomes empty or jug $C$ becomes full. Depending on the reason for the termination of pouring, we have the following subcases.
		
		\begin{subcase}
			Pouring terminated because the jug $A$ became empty.
		\end{subcase}
		
		In this case, since jug $A$ is empty, $d$ gallons of wine is distributed only in jugs $B$ and $C$ which implies that $i'+j'=d$. This together with the fact that $i'=i$ implies that $j'=d-i$. Thus we have $j' \in \{0,c,d-i,d-a-i\}$.
		
		\begin{subcase}
			Pouring terminated because the jug $C$ became full.
		\end{subcase}
		
		In this case, since the jug $C$ is full, the amount of wine in the jug $C$ is equal to its capacity. Hence, we have that $j'=c$ indicating that $j' \in \{0,c,d-i,d-a-i\}$.
		
		In any subcase, we have $j' \in \{0,c,d-i,d-a-i\}$, collectively implying $(\star)$, as desired.
		
		\begin{case}
			The measurable pour is $C \downarrow A$.
		\end{case}
		
		In this case, since the jug $B$ is untouched, we have $i'=i$. Since a measurable pour is being performed, we stop pouring only if jug $C$ becomes empty or jug $A$ becomes full. Again, depending on the reason for the termination of pouring, we have the following subcases.
		
		\begin{subcase}
			Pouring terminated because the jug $C$ became empty.
		\end{subcase}
		
		In this case, since the jug $C$ is empty, the amount of wine in the jug $C$ is equal to $0$. Hence, we have that $j'=0$ indicating that $j' \in \{0,c,d-i,d-a-i\}$.
		
		\begin{subcase}
			Pouring terminated because the jug $A$ became full.
		\end{subcase}
		
		In this case, since jug $A$ is full, the amount of wine in jug $A$ is equal to its capacity and the remaining wine is distributed over jugs $B$ and $C$ which implies that $a+i'+j'=d$. Further, substituting $i'=i$, we get $a+i+j'=d$. Hence, we have that $j'=d-a-i$ indicating that $j' \in \{0,c,d-i,d-a-i\}$.
		
		In any subcase, we have $j' \in \{0,c,d-i,d-a-i\}$, collectively implying $(\star)$, as desired.
		
		\begin{case}
			The measurable pour is $B \downarrow C$.
		\end{case}
		
		In this case, since the jug $A$ is untouched, we have $i+j=i'+j'$. Since a measurable pour is being performed, we stop pouring only if the jug $B$ becomes empty or the jug $C$ becomes full. Depending on the reason for the termination of pouring, we have the following subcases.
		
		\begin{subcase}
			Pouring terminated because the jug $B$ became empty.
		\end{subcase}
		
		In this case, since the jug $B$ is empty, the amount of wine in the jug $B$ is equal to $0$. Hence, we have that $i'=0$ indicating that $i' \in \{0,b\}$.
		
		\begin{subcase}
			Pouring terminated because the jug $C$ became full.
		\end{subcase}
		
		In this case, since the jug $C$ is full, the amount of wine in the jug $C$ is equal to its capacity. Hence, we have that $j'=c$ indicating that $j' \in \{0,c\}$.
		
		In any subcase, we have $i' \in \{0,b\}$ or $j' \in \{0,c\}$, collectively implying $(\star)$, as desired.
		
		\begin{case}
			The measurable pour is $C \downarrow B$.
		\end{case}
		
		In this case, since the jug $A$ is untouched, we have $i+j=i'+j'$. Since a measurable pour is being performed, we stop pouring only if jug $C$ becomes empty or jug $B$ becomes full. Depending on the reason for the termination of pouring, we have the following subcases.
		
		\begin{subcase}
			Pouring terminated because the jug $C$ became empty.
		\end{subcase}
		
		In this case, since the jug $C$ is empty, the amount of wine in the jug $C$ is equal to $0$. Hence, we have that $j'=0$ indicating that $j' \in \{0,c\}$.
		
		\begin{subcase}
			Pouring terminated because the jug $B$ became full.
		\end{subcase}
		
		In this case, since jug $B$ is full, the amount of wine in jug $B$ is equal to its capacity. Hence, we have that $i'=b$ indicating that $i' \in \{0,b\}$.
		
		In any subcase, we have $i' \in \{0,b\}$ or $j' \in \{0,c\}$, collectively implying $(\star)$, as desired.
		
		Since the cases are exhaustive and all the cases lead to the desired result, we can conclude that there exists a directed edge from the vertex $(i,j)$ to the vertex $(i',j')$ in $G_P$ which completes one direction of the proof.
		
		Conversely, assume that there exists a directed edge from a vertex $(i,j)$ to a vertex $(i',j')$ in $G_P$. Now, the goal is to prove that:
		\begin{equation}\tag{$\dagger$}
			\begin{split}
				& \text{the distribution } (i',j') \text{ can be obtained from the distribution } (i,j) \text{ in the puzzle } \mathbb{P}^{\tilde{a},\tilde{b},\tilde{c}}_{a,b,c}\\ & \text{by a single valid pour of wine from one jug to another.}
			\end{split}
		\end{equation}
		Since $G_P$ is a graph obtained using the model in Definition~\ref{def:Model} and since there is an edge from $(i,j)$ to $(i',j')$, the following conditions must hold:
		\begin{enumerate}[(1).]
			\item $0 \le d-i-j \le a$
			\item $0 \le d-i'-j' \le a$
		\end{enumerate}
		
		Further, the edge is present because one of the three conditions in the model holds. Suppose the first condition $(i)$ holds, then we have the following:
		\begin{enumerate}[(1).]
			\setItemnumber{3}
			\item $j=j'$
			\item $i' \in \{0,b,d-j,d-a-j\}$
		\end{enumerate}
		
		Since $(i,j)$ and $(i',j')$ are two different vertices or distributions, $(3)$ implies that $i \neq i'$. Note that the distribution $(i,j)$ implies that jugs $A$, $B$, and $C$ are filled with $d-i-j$, $i$, and $j$ gallons of wine respectively. Similarly, the distribution $(i',j')$ implies that jugs $A$, $B$, and $C$ are filled with $d-i'-j'$, $i'$, and $j'$ gallons of wine respectively.
		
		Suppose $i'=0$. Notice that to achieve the distribution $(0,j')$ from $(i,j)$ with $j'=j$, we need to pour the entire wine in jug $B$ to $A$. Since $i'=0$, we have $d-j' \le a$ by $(2)$. But combining this with $(3)$, we obtain $d-j \le a$, which indicates that in the distribution $(i,j)$, jug $A$ has sufficient space for the entire wine in jug $B$. Therefore, one can obtain the distribution $(i',j')$ from the distribution $(i,j)$ by pouring the wine from jug $B$ to jug $A$ until jug $B$ is empty, which results in a measurable pour $B \downarrow A$.
		
		Otherwise, suppose $i'=b$. Notice that to achieve the distribution $(b,j')$ from $(i,j)$ with $j'=j$, we need to pour the wine from jug $A$ to $B$ until $B$ is full. Observe that $(2)$ also implies that $i'+j' \le d$. This combined with the fact that $i'=b$ implies that $b+j' \le d$. Further, using $(3)$, we obtain $b+j \le d$ which in turn leads to $d-i-j \ge b-i$, indicating that in the distribution $(i,j)$, jug $A$ has sufficient wine to fill jug $B$ to its capacity, i.e., $b$. Therefore, one can obtain the distribution $(i',j')$ from the distribution $(i,j)$ by pouring the wine from jug $A$ to jug $B$ until jug $B$ is full, which results in a measurable pour $A \downarrow B$.
		
		Otherwise, suppose $i'=d-j$. Notice that to achieve the distribution $(d-j,j')$ from $(i,j)$ with $j'=j$, we need to pour the entire wine from jug $A$ to jug $B$. Since $i'$ was picked to be a non-negative integer inclusively between $0$ and $b$, we have $i' = d-j \le b$, indicating that in the distribution $(i,j)$, jug $B$ has sufficient space for the entire wine in jug $A$. Therefore, one can obtain the distribution $(i',j')$ from the distribution $(i,j)$ by pouring the wine from jug $A$ to jug $B$ until jug $A$ is empty, which results in a measurable pour $A \downarrow B$.
		
		Otherwise, by $(4)$ we have that $i'=d-a-j$. Notice that to achieve the distribution $(d-a-j,j')$ from $(i,j)$ with $j'=j$, we need to pour the wine from jug $B$ to jug $A$ until $A$ is full. From $i'=d-a-j$, we can deduce the following inequalities:
		\begin{align*}
			& i' = d-a-j \\
			\implies & d-j = a+i' \\
			\implies & a \le d-j \\
			\implies & a+i \le d+i-j \\
			\implies & i \ge a-d+i+j
		\end{align*}
		Thus we have $i \ge a-(d-i-j)$, indicating that in the distribution $(i,j)$, jug $B$ has sufficient wine to fill jug $A$ to its capacity, i.e., $a$, since the initial quantity of wine in jug $A$ is $d-i-j$. Therefore, one can obtain the distribution $(i',j')$ from the distribution $(i,j)$ by pouring the wine from jug $B$ to jug $A$ until jug $A$ is full, which results in a measurable pour $B \downarrow A$.
		
		Thus in any case, if the first condition $(i)$ holds, then $(\dagger)$ is true.
		
		Otherwise, suppose the second condition $(ii)$ holds, then we have the following:
		\begin{enumerate}[(1).]
			\setItemnumber{3}
			\item $i=i'$
			\item $j' \in \{0,c,d-i,d-a-i\}$
		\end{enumerate}
		
		Since $(i,j)$ and $(i',j')$ are two different vertices or distributions, $(3)$ implies that $j \neq j'$. Suppose $j'=0$. Notice that to achieve the distribution $(i',0)$ from $(i,j)$ with $i'=i$, we need to pour the entire wine in jug $C$ to $A$. Since $j'=0$, we have $d-i' \le a$ by $(2)$. But combining this with $(3)$, we obtain $d-i \le a$, which indicates that in the distribution $(i,j)$, jug $A$ has sufficient space for the entire wine in jug $C$. Therefore, one can obtain the distribution $(i',j')$ from the distribution $(i,j)$ by pouring the wine from jug $C$ to jug $A$ until jug $C$ is empty, which results in a measurable pour $C \downarrow A$.
		
		Otherwise, suppose $j'=c$. Notice that to achieve the distribution $(i',c)$ from $(i,j)$ with $i'=i$, we need to pour the wine from jug $A$ to $C$ until $C$ is full. Observe that $(2)$ also implies $i'+j' \le d$. This combined with the fact that $j'=c$, implies that $i'+c \le d$. Further, using $(3)$, we obtain $i+c \le d$ which in turn leads to $d-i-j \ge c-j$, indicating that in the distribution $(i,j)$, the jug $A$ has sufficient wine to fill the jug $C$ to its capacity, i.e., $c$. Therefore, one can obtain the distribution $(i',j')$ from the distribution $(i,j)$ by pouring the wine from jug $A$ to jug $C$ until the jug $C$ is full, which results in a measurable pour $A \downarrow C$.
		
		Otherwise, suppose $j'=d-i$. Notice that to achieve the distribution $(i',d-i)$ from $(i,j)$ with $i'=i$, we need to pour the entire wine from jug $A$ to jug $C$. Since $j'$ was picked to be a non-negative integer inclusively between $0$ and $c$, we have $j' = d-i \le c$, indicating that in the distribution $(i,j)$, jug $C$ has sufficient space for the entire wine in jug $A$. Therefore, one can obtain the distribution $(i',j')$ from the distribution $(i,j)$ by pouring the wine from jug $A$ to jug $C$ until jug $A$ is empty, which results in a measurable pour $A \downarrow C$.
		
		Otherwise, by $(4)$ we have that $j'=d-a-i$. Notice that to achieve the distribution $(i',d-a-i)$ from $(i,j)$ with $i'=i$, we need to pour the wine from jug $C$ to $A$ until $A$ is full. From $j'=d-a-i$, we can deduce the following inequalities:
		\begin{align*}
			& j'=d-a-i \\
			\implies & d-i = a+j' \\
			\implies & a \le d-i \\
			\implies & a+j \le d-i+j \\
			\implies & j \ge a-d+i+j
		\end{align*}
		Thus we have $j \ge a-(d-i-j)$, indicating that in the distribution $(i,j)$, jug $C$ has sufficient wine to fill jug $A$ to its capacity, i.e., $a$, since the initial quantity of wine in jug $A$ is $d-i-j$. Therefore, one can obtain the distribution $(i',j')$ from the distribution $(i,j)$ by pouring the wine from jug $C$ to jug $A$ until jug $A$ is full, which results in a measurable pour $C \downarrow A$.
		
		Thus in any case, if the second condition $(ii)$ holds, then $(\dagger)$ is true.
		
		Otherwise, if the conditions $(i)$ and $(ii)$ do not hold, then by Definition~\ref{def:Model}, the third condition $(iii)$ holds. Thus we have the following:
		\begin{enumerate}[(1).]
			\setItemnumber{3}
			\item $i+j=i'+j'$
			\item $i' \in \{0,b\}$ or $j' \in \{0,c\}$
		\end{enumerate}
		
		Suppose $i'=0$. Notice that to achieve the distribution $(0,j')$ from $(i,j)$ satisfying $(3)$, we need to pour the entire wine from jug $B$ to jug $C$. Since $i'=0$, we have $i+j=j'$ by $(3)$. Therefore, since $j'$ was picked to be a non-negative integer inclusively between $0$ and $c$, we have $c \ge j' = i+j$ which in turn implies that $c-j \ge i$, indicating that in the distribution $(i,j)$, the jug $C$ has sufficient space for the entire wine in the jug $B$. Therefore, one can obtain the distribution $(i',j')$ from the distribution $(i,j)$ by pouring the wine from jug $B$ to jug $C$ until jug $B$ is empty, which results in a measurable pour $B \downarrow C$.
		
		Otherwise, suppose $i'=b$. Notice that to achieve the distribution $(b,j')$ from $(i,j)$ satisfying $(3)$, we need to pour the wine from jug $C$ to $B$ until $B$ is full. The fact that $i'=b$ together with $(3)$ implies that $b = i' \le i'+j' = i+j$. Thus we have $j \ge b-i$, indicating that in the distribution $(i,j)$, the jug $C$ has sufficient wine to fill the jug $B$ to its capacity, i.e., $b$. Therefore, one can obtain the distribution $(i',j')$ from the distribution $(i,j)$ by pouring the wine from jug $C$ to jug $B$ until jug $B$ is full, which results in a measurable pour $C \downarrow B$.
		
		Otherwise, suppose $j'=0$. Notice that to achieve the distribution $(i',0)$ from $(i,j)$ satisfying $(3)$, we need to pour the entire wine from jug $C$ to jug $B$. Since $j'=0$, we have $i+j=i'$ by $(3)$. Therefore, since $i'$ was picked to be a non-negative integer inclusively between $0$ and $b$, we have $b \ge i' = i+j$ which in turn implies that $b-i \ge j$ indicates that in the distribution $(i,j)$, the jug $B$ has sufficient space for the entire wine in the jug $C$. Therefore, one can obtain the distribution $(i',j')$ from the distribution $(i,j)$ by pouring the wine from jug $C$ to jug $B$ until jug $C$ is empty, which results in a measurable pour $C \downarrow B$.
		
		Otherwise, if $i' \notin \{0,b\}$ and $j' \neq 0$, then by $(4)$, we are sure that $j'=c$. Notice that to achieve the distribution $(i',c)$ from $(i,j)$ satisfying $(3)$, we need to pour the wine from jug $B$ to $C$ until $C$ is full. The fact that $j'=c$ together with $(3)$ implies that $c = j' \le i'+j' = i+j$. Thus we have $i \ge c-j$, indicating that in the distribution $(i,j)$, the jug $B$ has sufficient wine to fill the jug $C$ to its capacity, i.e., $c$. Therefore, one can obtain the distribution $(i',j')$ from the distribution $(i,j)$ by pouring the wine from jug $B$ to jug $C$ until jug $C$ is full, which results in a measurable pour $B \downarrow C$.
		
		Hence, collectively, if any one of the three conditions in the model holds, then $(\dagger)$ is true. Since there is a directed edge from the vertex $(i,j)$ to the vertex $(i',j')$, one of the three conditions in the model holds by Definition~\ref{def:Model}, which marks the completion of the proof.
	\end{proof}
	
	Once we have Theorem~\ref{thm:ThreeJug}, the following corollary follows from the same. This corollary gives an equivalent condition for the three-jug puzzle to have a solution, by using the model in Definition~\ref{def:Model}.
	
	\begin{corollary}\label{cor:ThreeJug}
		Let $\mathbb{P}^{\tilde{a},\tilde{b},\tilde{c}}_{a,b,c}$ be a generic three-jug puzzle and let $G_P$ be the graph obtained from the ordered quadruple $P=(a,b,c,\tilde{a}+\tilde{b}+\tilde{c})$ using the graph model in Definition~\ref{def:Model}. Then the puzzle $\mathbb{P}^{\tilde{a},\tilde{b},\tilde{c}}_{a,b,c}$ has a solution if and only if there exists a directed path from the vertex $(\tilde{b},\tilde{c})$ to the vertex $(\frac{\tilde{a}+\tilde{b}+\tilde{c}}{2},0)$ in $G_P$.
	\end{corollary}
	
	\begin{proof}
		First, assume that the puzzle $\mathbb{P}^{\tilde{a},\tilde{b},\tilde{c}}_{a,b,c}$ has a solution. Then, one can reach the distribution $(\frac{\tilde{a}+\tilde{b}+\tilde{c}}{2},0)$ from the distribution $(\tilde{b},\tilde{c})$ through a sequence of intermediate distributions obtained exactly by the sequence of measurable pours performed along the process. For every single change of distribution in this sequence of pouring, we have a directed edge between the corresponding vertices, by Theorem~\ref{thm:ThreeJug}. By taking the union of all these directed edges corresponding to the entire sequence that is used to reach the distribution $(\frac{\tilde{a}+\tilde{b}+\tilde{c}}{2},0)$ from the distribution $(\tilde{b},\tilde{c})$, we can obtain a directed path from the vertex $(\tilde{b},\tilde{c})$ to the vertex $(\frac{\tilde{a}+\tilde{b}+\tilde{c}}{2},0)$ in $G_P$, as desired. This completes one direction of the proof.
		
		Conversely, assume that there exists a directed path from the vertex $(\tilde{b},\tilde{c})$ to the vertex $(\frac{\tilde{a}+\tilde{b}+\tilde{c}}{2},0)$ in $G_P$. But a directed path is a sequence of directed edges. For each edge in this sequence, by Theorem~\ref{thm:ThreeJug}, we can obtain the corresponding ordered pair of distributions such that the latter can be obtained from the first by a single measurable pour of the wine from one jug to another. This sequence of ordered pairs starting from $(\tilde{b},\tilde{c})$ ending at $(\frac{\tilde{a}+\tilde{b}+\tilde{c}}{2},0)$, generates the sequential order of pouring of wine between the jugs which gives us the desired solution to the puzzle $\mathbb{P}^{\tilde{a},\tilde{b},\tilde{c}}_{a,b,c}$. This completes the proof of the corollary.
	\end{proof}
	
	\section{Sketch of an Algorithm}
	
	Now, if we are provided with a generic three-jug puzzle, one can just extract an ordered quadruple out of it and obtain the graph model presented in Definition~\ref{def:Model}. Once that is done, it all boils down to determining whether the puzzle has a solution or equivalently whether the graph has the corresponding directed path. In other words, by Corollary~\ref{cor:ThreeJug}, we can leave aside the puzzle and just look at the fairly simple graph model and look at its structural properties and infer whether the original generic three-jug puzzle had a solution or not; if it had a solution, then the same can be determined using the graph model.
	
	Once we are sure that the puzzle has a solution, then one can think of the problem of minimizing the number of measurable pours, i.e., to determine the minimum number of measurable pours that is necessary to obtain a solution. Again, by using Corollary~\ref{cor:ThreeJug}, this boils down to determining the shortest directed path in the graph model. Observe that whenever a directed path exists in a graph, a shortest directed path also exists in the graph. Alternatively, whenever a solution exists for the puzzle, a solution with the minimum number of measurable pours also exists for the puzzle. Therefore, the model can also be used to obtain the optimum solution, if it exists.
	
	As a result of the previous discussions, given a generic three-jug puzzle $\mathbb{P}^{\tilde{a},\tilde{b},\tilde{c}}_{a,b,c}$, one can use the following steps to obtain a solution. These steps can be treated as the skeleton of an algorithm.
	
	\begin{enumerate}[1.]
		\item Obtain an ordered quadruple $Q=(a,b,c,\tilde{a}+\tilde{b}+\tilde{c})$ from the puzzle $\mathbb{P}^{\tilde{a},\tilde{b},\tilde{c}}_{a,b,c}$.
		
		\item Obtain a directed graph $G_Q$ corresponding to $Q$ using the model in Definition~\ref{def:Model}.
		
		\item Look for a directed path from the vertex $s = (\tilde{b},\tilde{c})$ to the vertex $t = (\frac{\tilde{a}+\tilde{b}+\tilde{c}}{2},0)$ in $G_Q$.
		
		\item If there exists a directed $st$-path in $G_Q$, then the vertex sequence in the same order as a shortest directed $st$-path will be a solution to the puzzle $\mathbb{P}^{\tilde{a},\tilde{b},\tilde{c}}_{a,b,c}$.
		
		\item Otherwise, if there is no directed $st$-path in $G_Q$, conclude that the puzzle $\mathbb{P}^{\tilde{a},\tilde{b},\tilde{c}}_{a,b,c}$ has no solution.
	\end{enumerate}
	
	\section{Conclusion}
	
	The Three-Jug Puzzle, once a recreational curiosity, has now evolved into a rich domain for exploring algorithm design, cognitive reasoning, algebraic modeling, and combinatorial logic. Starting from heuristic solutions and geometric visualizations to formal complexity analysis and computational approaches, there is a huge evolution in the literature. Even, determining whether these kinds of puzzles have a solution or not is an interesting research problem to be tackled. Several variants of the puzzle have been shown to be solvable using various methods and algorithms. In view of this paper, a generic version of the three-jug puzzle is considered and an independent directed graph model is proposed such that the solvability of the puzzle can be determined by using the model. If a solution exists, then the same can be determined by using the proposed model. It is important to note that the model as such does not depend on the puzzle.
	
	In the future, one can consider a different variant of the puzzle or a completely different puzzle and try to determine the solvability of the same by using a similar technique or the idea of the model presented in the paper. Our motive was to utilize graph theory as a tool to tackle the generic three-jug puzzle. In the future, researchers may utilize some other field of mathematics to tackle this puzzle and try to come up with a better model, which leaves a wide scope ahead.
	
	\bibliographystyle{SK}
	\bibliography{ThreeJug}
	
\end{document}